\documentclass[10pt]{amsart}
\usepackage{amsfonts}
\usepackage{ifthen}
\usepackage{amsthm}
\usepackage{amsmath}
\usepackage{graphicx}
\usepackage{amscd,amssymb,amsthm}

\newcounter{minutes}
\divide\time by 60
\newcounter{hours}
\multiply\time by 60 \addtocounter{minutes}{-\time}

\newcommand{\real}{\operatorname{Re}}

\newtheorem{theorem}{Theorem}


\title[\'A. Baricz, D. Jankov, T.K. Pog\'any, R. Sz\'asz/Identity for zeros of Bessel functions]{On an identity for zeros of Bessel functions}

\author[]{\'Arp\'ad Baricz}
\address{Department of Economics, Babe\c{s}-Bolyai University, 400591 Cluj-Napoca, Romania and Institute of Applied Mathematics, John von Neumann Faculty of Informatics, \'Obuda University, 1034 Budapest, Hungary} \email{bariczocsi@yahoo.com}

\author[]{Dragana Jankov Ma\v sirevi\' c}
\address{Department of Mathematics, University of Osijek, 31000 Osijek, Croatia} \email{djankov@mathos.hr}

\author[]{Tibor K. Pog\'any}
\address{Faculty of Maritime Studies, University of Rijeka, 51000 Rijeka, Croatia} \email{poganj@pfri.hr}

\author[]{R\'obert Sz\'asz}
\address{Department of Mathematics and Informatics, Sapientia Hungarian University of Transylvania, 540485 T\^argu-Mure\c{s}, Romania}
\email{rszasz@ms.sapientia.ro}

\begin{document}

\def\thefootnote{}
\footnotetext{ \texttt{File:~\jobname .tex,
          printed: \number\year-0\number\month-0\number\day,
          \thehours.\ifnum\theminutes<10{0}\fi\theminutes}
} \makeatletter\def\thefootnote{\@arabic\c@footnote}\makeatother

\maketitle

\begin{abstract}
In this paper our aim is to present an elementary proof of an identity of Calogero concerning the zeros of Bessel functions of the first kind. Moreover, by using our elementary approach we present a new identity for the zeros of Bessel functions of the first kind, which in particular reduces to some other new identities. We also show that our method can be applied for the zeros of other special functions, like Struve functions of the first kind, and modified Bessel functions of the second kind.
\end{abstract}

\section{Introduction and Main Results}

In 1977 F. Calogero \cite{calogero} deduced the following identity
\begin{equation}\label{formCalo}
\sum_{n\geq1,n\neq{k}}\frac{1}{j_{\nu,n}^2-{j_{\nu,k}^2}}=\frac{\nu+1}{2j_{\nu,k}^2},
\end{equation}
where $\nu>-1,$ $k\in\{1,2,\dots\}$ and $j_{\nu,n}$ stands for the $n$th positive zero of the Bessel function of the first kind $J_{\nu}.$ Calogero's proof \cite{calogero} of \eqref{formCalo} is based on the infinite product representation of the Bessel functions of the first kind and on the clever use of an equivalent form of the Mittag-Leffler expansion
\begin{equation}\label{mittag}
\frac{J_{\nu+1}(x)}{J_{\nu}(x)}=\sum_{n\geq 1}\frac{2x}{j_{\nu,n}^2-x^2},
\end{equation}
where $\nu>-1.$ Note that in \cite{calogero} it was pointed out that results like \eqref{formCalo} are related to
the connection between the motion of poles and zeros of special solutions of
partial differential equations and many-body problems. In 1986 Ismail and Muldoon \cite{ismail} mentioned that \eqref{formCalo} can be obtained also by evaluating the residues of the functions in \eqref{formCalo} at their poles. In this paper our aim is to present an alternative proof of \eqref{formCalo} by using only elementary analysis. Our proof is based on the Mittag-Leffler expansion \eqref{mittag}, recurrence relations, the Bessel differential equation and on the Bernoulli-Hospital rule for the limit of quotients. Moreover, by using our idea we are able to prove the following new results.

\begin{theorem}\label{thm1}
If $\nu>-1$ and $k\in\{1,2,\dots\},$ then we have
\begin{equation}\label{formnew}
\sum_{\substack{{n\geq1},{n\neq{k}}}}\frac{1}{j_{\nu,n}^4-{j_{\nu,k}^4}}=-\frac{1}{2j_{\nu,k}^2}\sum_{n\geq 1}\frac{1}{j_{\nu,n}^2+j_{\nu,k}^2}+\frac{\nu+2}{4j_{\nu,k}^4}.
\end{equation}
In particular, for all $k\in\{1,2,\dots\}$ we have
\begin{equation}\label{formnew1/2}
\sum_{\substack{{n\geq1},{n\neq{k}}}}\frac{1}{{n}^4-{{k}^4}}=-\frac{\pi}{4k^3}\coth(k\pi)+\frac{7}{8k^4}.
\end{equation}
Moreover, for all $k\in\{1,3,\dots\}$ we have
\begin{equation}\label{formnew-1/2}
\sum_{\substack{{n\geq1},{n\neq{k}}\\ n\ is\ odd}}\frac{1}{{n}^4-{{k}^4}}=-\frac{\pi}{8k^3}\tanh\left(\frac{k\pi}{2}\right)+\frac{3}{8k^4}.
\end{equation}
\end{theorem}

As far as we know the above results are new and as we can see below our method can be applied for the zeros of other special functions, like Struve functions and modified Bessel functions of the second kind. Our proof in this case is based on the corresponding Mittag-Leffler expansion for Struve functions, recurrence relations, the Struve differential equation and on the Bernoulli-Hospital rule. During the process of writing this paper an alternative proof of \eqref{formnew} was proposed by Christophe Vignat. His proof was based on the formula
\begin{equation}\label{vignat}
\sum_{n\geq 1,n\neq k}\frac{1}{j_{\nu,n}^4-j_{\nu,k}^4}=\frac{1}{2j_{\nu,k}^2}\sum_{n\geq 1,n\neq k}\frac{1}{j_{\nu,n}^2-j_{\nu,k}^2}-\frac{1}{2j_{\nu,k}^2}\sum_{n\geq 1,n\neq k}\frac{1}{j_{\nu,n}^2+j_{\nu,k}^2},
\end{equation}
which by means of \eqref{formCalo} evidently implies \eqref{formnew}. It is worth to mention that \eqref{vignat} holds true in fact for any set of numbers $\{j_{\nu,n}\}_{n\geq1},$ and by using this idea we can get \eqref{formStruve4} for the zeros of Struve functions.

\begin{theorem}\label{thm2}
Let $h_{\nu,n}$ be the $n$th positive zero of the Struve function of the first kind $\mathbf{H}_{\nu}.$
If $|\nu|<\frac{1}{2}$ and $k\in\{1,2,\dots\},$ then the following identities are valid
\begin{equation}\label{formStruve}
\sum_{{{n\geq1},{n\neq{k}}}}\frac{1}{h_{\nu,n}^2-{h_{\nu,k}^2}}=
\frac{\nu+2}{2h_{\nu,k}^2}-\frac{h_{\nu,k}^{\nu-2}}{\sqrt{\pi}2^{\nu-1}\Gamma\left(\nu+\frac{1}{2}\right)\mathbf{H}_{\nu}'(h_{\nu,k})},
\end{equation}
\begin{equation}\label{formStruve4}
\sum_{{{n\geq1},{n\neq{k}}}}\frac{1}{h_{\nu,n}^4-{h_{\nu,k}^4}}=
-\frac{1}{2h_{\nu,k}^2}\sum_{{{n\geq1}}}\frac{1}{h_{\nu,n}^2+{h_{\nu,k}^2}}+
\frac{\nu+3}{4h_{\nu,k}^4}-\frac{h_{\nu,k}^{\nu-4}}{\sqrt{\pi}2^{\nu}\Gamma\left(\nu+\frac{1}{2}\right)\mathbf{H}_{\nu}'(h_{\nu,k})}.
\end{equation}
\end{theorem}

It is also worth to mention that since \cite[p. 291]{nist}
$\mathbf{H}_{-\frac{1}{2}}(x)=\sqrt{\frac{2}{\pi x}}\sin x,$
by tending with $\nu$ to $-\frac{1}{2}$ in \eqref{formStruve} and taking into account that $h_{-\frac{1}{2},n}=n\pi$ we obtain for all $k\in\{1,2,\dots\}$ the
following known result
\begin{equation}\label{formSnew-1/2}
\sum_{{{n\geq1},{n\neq{k}}}}\frac{1}{{n}^2-{{k}^2}}=\frac{3}{4k^2}.
\end{equation}
Note that this result can be obtained also from \eqref{formCalo} by choosing $\nu=\frac{1}{2},$ as it was pointed out in \cite{calogero}.

The properties of the zeros of modified Bessel function of the
second kind $K_{\nu}$, also called Macdonald function, were studied
by several authors, but not so often and in detail as the zeros of
the Bessel functions of the first and second kind, i.e. $J_{\nu}$
and $Y_{\nu}.$ Macdonald was the first who discussed the zeros of
that function and he showed \cite[p. 511]{Watson} that $K_\nu(z)$
has no positive zeros when $\nu\ge 0$ and also that it has no
zeros for which $|{\rm{arg}}z|\le\frac\pi2$. Note also that Ferreira
and Sesma \cite{FS} studied the zeros of $K_{\nu}$ when
$\nu\in\mathbb{C}$. In the sequel we are interested on the zeros of
$K_{\nu}$, when $\nu=n+\frac12$ and $n\in\mathbb{N}$. We know that (see for example \cite{YH}) $K_{n+\frac{1}{2}}$ for $n\in\{1,2,\dots\}$ has exactly $n$ zeros with
multiplicity one in $\mathbb{C}^-:=\{z\in\mathbb{C}:\real z<0\}$
and that the non--real zeros are complex conjugate in pairs
(see \cite{HMS}). We also know that for $\nu=n+\frac12$,
$n\in\mathbb{N}$ the value $K_\nu(z)$ can be represented as
$$
K_\nu(z)=\sqrt{\dfrac{\pi}{2z}}\dfrac{{\rm e}^{-z}}{z^n}H_n(z),
$$
where $H_n$ is a monic polynomial of degree $n$ defined as
$$H_n(z)=\sum_{k=0}^n\dfrac{\Gamma(2n-k+1)\,z^k}{\Gamma(n-k+1)\,\Gamma(k+1)\,2^{n-k}}\,.$$
If we denote the zeros of $K_\nu$ with
$z_{\nu,1},\dots,z_{\nu,n}$, observing that those are also the
zeros of $H_n$ we can conclude that (see \cite{YH})
$$H_n(z)=\prod_{j=1}^n(z-z_{\nu,j}).$$
By using the ideas of the proofs of the above theorems we are able to prove the next results.

\begin{theorem}\label{th3}
For all $\nu=n+\frac12$, $n\in\mathbb{N}$ and $j\in\{1,\dots,n\}$ we have
\begin{equation}\label{eqk1}
 \sum_{{k=1, k\neq j}}^n\dfrac{1}{z_{\nu,k}-z_{\nu,j}}=\dfrac{1-2z_{\nu,j}-2\nu}{2z_{\nu,j}},
\end{equation}
\begin{equation}\label{eqk2}
\sum_{{k=1,k\neq j}}^n\dfrac{1}{z_{\nu,k}^2-z_{\nu,j}^2}=\dfrac{1-z_{\nu,j}-\nu}{2z_{\nu,j}^2}-\frac{1}{2z_{\nu,j}}\sum_{k=1}^n\dfrac{1}{z_{\nu,k}+z_{\nu,j}},
\end{equation}
\begin{equation}\label{eqk4}
 \sum_{{k=1,k\neq j}}^n\dfrac{1}{z_{\nu,k}^4-z_{\nu,j}^4}=\dfrac{2-\nu-z_{\nu,j}}{4z_{\nu,j}^4}-\dfrac{1}{4z_{\nu,j}^3}
 \sum_{k=1}^n\left(\dfrac{1}{z_{\nu,j}+z_{\nu,k}}+\dfrac{2z_{\nu,j}}{z_{\nu,j}^2+z_{\nu,k}^2}\right).
\end{equation}
\end{theorem}

\section{Proofs}
\setcounter{equation}{0}

In this section we are going to present the proof of \eqref{formCalo} and of Theorems \ref{thm1}, \ref{thm2} and \ref{th3}.

\begin{proof}[\bf Proof of \eqref{formCalo}] We start with the following identities
\begin{equation}\label{eq1}{j_{\nu,k}}J_\nu''(j_{\nu,k})+J_\nu'(j_{\nu,k})=0,\end{equation}
\begin{equation}\label{eq2}{J}_{\nu+1}(j_{\nu,k})+{J}_\nu'(j_{\nu,k})=0,\end{equation}
\begin{equation}\label{eq3}{j_{\nu,k}}{J}_{\nu+1}'(j_{\nu,k})-(\nu+1){J}_{\nu}'(j_{\nu,k})=0,\end{equation}
which readily follow from the fact that $J_{\nu}$ is a particular solution of the Bessel differential equation \cite[p. 217]{nist}, that is, satisfies
$$x^2J_\nu''(x)+xJ_\nu'(x)+(x^2-\nu^2)J_\nu(x)=0,$$
and from the recurrence relations \cite[p. 222]{nist}
\begin{equation}\label{rec1}xJ_\nu'(x)-\nu{J}_\nu(x)=-xJ_{\nu+1}(x),\end{equation}
$$xJ_{\nu+1}'(x)+(\nu+1){J}_{\nu+1}(x)=xJ_{\nu}(x).$$
By using the Mittag-Leffler expansion \eqref{mittag} we obtain that
$$\mathbf{\Omega}_1=\sum_{\substack{{n\geq1}\\{n\neq{k}}}}\frac{1}{j_{\nu,n}^2-{j_{\nu,k}^2}}
=\lim_{x\rightarrow{j_{\nu,k}}}\left(\frac{J_{\nu+1}(x)}{2xJ_{\nu}(x)}-\frac{1}{j_{\nu,k}^2-x^2}\right)=
\frac{1}{2j_{\nu,k}}\lim_{x\rightarrow{j_{\nu,k}}}\frac{(j_{\nu,k}^2-x^2)J_{\nu+1}(x)-2xJ_{\nu}(x)}{(j_{\nu,k}^2-x^2)J_{\nu}(x)}.$$
Now, applying the Bernoulli-Hospital rule two times and the relations \eqref{eq1}, \eqref{eq2} and \eqref{eq3} we get
\begin{align*}\mathbf{\Omega}_1&=\frac{1}{2j_{\nu,k}}\lim_{x\rightarrow{j_{\nu,k}}}\frac{(j_{\nu,k}^2-x^2)J_{\nu+1}''(x)-4xJ_{\nu+1}'(x)-2J_{\nu+1}(x)-4J_{\nu}'(x)-2xJ_{\nu}''(x)}
{(j_{\nu,k}^2-x^2)J_{\nu}''(x)-4xJ_{\nu}'(x)-2J_{\nu}(x)}\\
&=\frac{1}{2j_{\nu,k}}\left(\frac{-4j_{\nu,k}J_{\nu+1}'(j_{\nu,k})-2J_{\nu+1}(j_{\nu,k})-4J_{\nu}'(j_{\nu,k})-2j_{\nu,k}J_{\nu}''(j_{\nu,k})}
{-4j_{\nu,k}J_{\nu}'(j_{\nu,k})}\right)=\frac{\nu+1}{2j_{\nu,k}^2}.
\end{align*}
\end{proof}

\begin{proof}[\bf Proof of Theorem \ref{thm1}] Let us denote with $I_{\nu}$ the modified Bessel function of the first kind or Bessel function of the first kind with purely imaginary argument. By using the Weierstrassian products \cite[p. 235]{nist}
$$2^{\nu}\Gamma(\nu+1)x^{-\nu}J_{\nu}(x)=\prod_{n\geq 1}\left(1-\frac{x^2}{j_{\nu,n}^2}\right), \ \
2^{\nu}\Gamma(\nu+1)x^{-\nu}I_{\nu}(x)=\prod_{n\geq 1}\left(1+\frac{x^2}{j_{\nu,n}^2}\right),$$
we obtain
$$2^{2\nu}\Gamma^2(\nu+1)x^{-2\nu}J_{\nu}(x)I_{\nu}(x)=\prod_{n\geq 1}\left(1-\frac{x^4}{j_{\nu,n}^4}\right).$$
Logarithmic differentiation gives
$$-\frac{2\nu}{x}+\frac{J_{\nu}'(x)}{J_{\nu}(x)}+\frac{I_{\nu}'(x)}{I_{\nu}(x)}=-\sum_{n\geq1}\frac{4x^3}{j_{\nu,n}^4-x^4},$$
which in view of \eqref{rec1} and its analogue
$$xI_\nu'(x)-\nu{I}_\nu(x)=xI_{\nu+1}(x),$$
can be rewritten as
$$\frac{1}{4x^3}\left(\frac{J_{\nu+1}(x)}{J_{\nu}(x)}-\frac{I_{\nu+1}(x)}{I_{\nu}(x)}\right)=\sum_{n\geq1}\frac{1}{j_{\nu,n}^4-x^4}.$$
By using the above Mittag-Leffler expansion we obtain that
\begin{align*}\mathbf{\Omega}_2&=\sum_{\substack{{n\geq1}\\{n\neq{k}}}}\frac{1}{j_{\nu,n}^4-{j_{\nu,k}^4}}
=\lim_{x\rightarrow{j_{\nu,k}}}\left(\frac{J_{\nu+1}(x)}{4x^3J_{\nu}(x)}-\frac{I_{\nu+1}(x)}{4x^3I_{\nu}(x)}-\frac{1}{j_{\nu,k}^4-x^4}\right)\\
&=\lim_{x\rightarrow{j_{\nu,k}}}\left(\frac{J_{\nu+1}(x)}{4x^3J_{\nu}(x)}-\frac{1}{j_{\nu,k}^4-x^4}\right)-
\frac{I_{\nu+1}(j_{\nu,k})}{4j_{\nu,k}^3I_{\nu}(j_{\nu,k})}\\
&=\lim_{x\rightarrow{j_{\nu,k}}}\frac{1}{(j_{\nu,k}^2+x^2)x^2}\left(\frac{(j_{\nu,k}^2+x^2)J_{\nu+1}(x)}{4xJ_{\nu}(x)}-\frac{x^2}{j_{\nu,k}^2-x^2}\right)-
\frac{I_{\nu+1}(j_{\nu,k})}{4j_{\nu,k}^3I_{\nu}(j_{\nu,k})}.
\end{align*}
Now, applying again the Bernoulli-Hospital rule two times and the relations \eqref{eq1}, \eqref{eq2} and \eqref{eq3} we get
\begin{align*}\mathbf{\Omega}_3&=\lim_{x\rightarrow{j_{\nu,k}}}\frac{1}{(j_{\nu,k}^2+x^2)x^2}\left(\frac{(j_{\nu,k}^2+x^2)J_{\nu+1}(x)}{4xJ_{\nu}(x)}-\frac{x^2}{j_{\nu,k}^2-x^2}\right)\\
&=\frac{1}{8j_{\nu,k}^5}\lim_{x\rightarrow{j_{\nu,k}}}\frac{(j_{\nu,k}^4-x^4)J_{\nu+1}(x)-4x^3J_{\nu}(x)}{(j_{\nu,k}^2-x^2)J_{\nu}(x)}\\
&=\frac{1}{8j_{\nu,k}^5}\lim_{x\rightarrow{j_{\nu,k}}}\frac{(j_{\nu,k}^4-x^4)J_{\nu+1}''(x)-8x^3J_{\nu+1}'(x)-12x^2J_{\nu+1}(x)-24x(J_{\nu}(x)+xJ_{\nu}'(x))-4x^3J_{\nu}''(x)}
{(j_{\nu,k}^2-x^2)J_{\nu}''(x)-4xJ_{\nu}'(x)-2J_{\nu}(x)}\\
&=\frac{1}{8j_{\nu,k}^5}\left(\frac{-8j_{\nu,k}^3J_{\nu+1}'(j_{\nu,k})-12j_{\nu,k}^2J_{\nu+1}(j_{\nu,k})-24j_{\nu,k}^2J_{\nu}'(j_{\nu,k})-4j_{\nu,k}^3J_{\nu}''(j_{\nu,k})}
{-4j_{\nu,k}J_{\nu}'(j_{\nu,k})}\right)=\frac{\nu+2}{4j_{\nu,k}^4}.
\end{align*}
Finally, by using the Mittag-Leffler expansion
\begin{equation}\label{leffler}\frac{I_{\nu+1}(x)}{I_{\nu}(x)}=\sum_{n\geq 1}\frac{2x}{j_{\nu,n}^2+x^2}\end{equation}
the proof of \eqref{formnew} is complete. Now, by taking $\nu=\frac{1}{2}$ and $\nu=-\frac{1}{2}$ in \eqref{formnew} we get \eqref{formnew1/2} and
$$\sum_{\substack{{n\geq1}\\{n\neq{k}}}}\frac{1}{{(2n-1)}^4-{{(2k-1)}^4}}=-\frac{\pi}{8(2k-1)^3}\tanh\left(\frac{(2k-1)\pi}{2}\right)+\frac{3}{8(2k-1)^4},$$
which is equivalent to \eqref{formnew-1/2}. Here we used the Mittag-Leffler expansion \cite[p. 126]{nist}
$$\coth (x)=\frac{1}{x}+2x\sum_{n\geq 1}\frac{1}{x^2+n^2\pi^2}$$
and \eqref{leffler} for $\nu=-\frac{1}{2}$ together with \cite[p. 254]{nist}
$$I_{\frac{1}{2}}(x)=\sqrt{\frac{2}{\pi x}}\sinh x,\ \ I_{-\frac{1}{2}}(x)=\sqrt{\frac{2}{\pi x}}\cosh x,\ \ \tanh x=\frac{I_{\frac{1}{2}}(x)}{I_{-\frac{1}{2}}(x)}.$$
Moreover, we have used the fact that for $n\in\{1,2,\dots\}$ we have $j_{\frac{1}{2},n}=n\pi$ and $j_{-\frac{1}{2},n}=\frac{(2n-1)\pi}{2}$ since $J_{\frac{1}{2}}(x)$ is proportional to $\sin x$ and $J_{-\frac{1}{2}}(x)$ is proportional to $\cos x,$ that is, we have \cite[p. 228]{nist}
$$J_{\frac{1}{2}}(x)=\sqrt{\frac{2}{\pi x}}\sin x,\ \ J_{-\frac{1}{2}}(x)=\sqrt{\frac{2}{\pi x}}\cos x.$$
\end{proof}

\begin{proof}[\bf Proof of Theorem \ref{thm2}] The proof is quite similar to the proof of \eqref{formCalo}. We start with the following identities
\begin{equation}\label{eq1h}{h_{\nu,k}}\mathbf{H}_\nu''(h_{\nu,k})+\mathbf{H}_\nu'(h_{\nu,k})=\frac{h_{\nu,k}^{\nu}}{\sqrt{\pi}2^{\nu-1}\Gamma\left(\nu+\frac{1}{2}\right)},\end{equation}
\begin{equation}\label{eq2h}\mathbf{H}_{\nu-1}(h_{\nu,k})-\mathbf{H}_\nu'(h_{\nu,k})=0,\end{equation}
\begin{equation}\label{eq3h}{h_{\nu,k}}\mathbf{H}_{\nu-1}'(h_{\nu,k})-\nu\mathbf{H}_{\nu}'(h_{\nu,k})=h_{\nu,k}\mathbf{H}_{\nu}''(h_{\nu,k}),\end{equation}
which readily follow from the fact that $\mathbf{H}_{\nu}$ is a particular solution of the Struve differential equation \cite[p. 288]{nist}, that is, satisfies
$$x\mathbf{H}_\nu''(x)+\mathbf{H}_\nu'(x)+x\left(1-\frac{\nu^2}{x^2}\right)\mathbf{H}_\nu(x)=\frac{x^{\nu}}{\sqrt{\pi}2^{\nu-1}\Gamma\left(\nu+\frac{1}{2}\right)},$$
and from the recurrence relations \cite[p. 292]{nist}
\begin{equation}\label{rec1h}x\mathbf{H}_\nu'(x)+\nu\mathbf{H}_\nu(x)=x\mathbf{H}_{\nu-1}(x),\end{equation}
$$\mathbf{H}_{\nu-1}(x)+x\mathbf{H}_{\nu-1}'(x)=(\nu+1)\mathbf{H}_{\nu}'(x)+x\mathbf{H}_{\nu}''(x).$$
By using the Mittag-Leffler expansion \cite[Lemma 1]{bps}
$$\frac{\mathbf{H}_{\nu-1}(x)}{\mathbf{H}_{\nu}(x)}=\frac{2\nu+1}{x}-\sum_{n\geq1}\frac{2x}{h_{\nu,n}^2-x^2}$$
we obtain that
\begin{align*}\mathbf{\Omega}_4&=\sum_{\substack{{n\geq1}\\{n\neq{k}}}}\frac{1}{h_{\nu,n}^2-{h_{\nu,k}^2}}
=\lim_{x\rightarrow{h_{\nu,k}}}\left(\frac{2\nu+1}{2x^2}-\frac{\mathbf{H}_{\nu-1}(x)}{2x\mathbf{H}_{\nu}(x)}-\frac{1}{h_{\nu,k}^2-x^2}\right)\\&=
\frac{2\nu+1}{2h_{\nu,k}^2}-\frac{1}{2h_{\nu,k}}\lim_{x\rightarrow{h_{\nu,k}}}
\frac{(h_{\nu,k}^2-x^2)\mathbf{H}_{\nu-1}(x)+2x\mathbf{H}_{\nu}(x)}{(h_{\nu,k}^2-x^2)\mathbf{H}_{\nu}(x)}.\end{align*}
Now, applying the Bernoulli-Hospital rule two times and the relations \eqref{eq1h}, \eqref{eq2h} and \eqref{eq3h} we get
\begin{align*}\mathbf{\Omega}_5&=\frac{1}{2h_{\nu,k}}\lim_{x\rightarrow{h_{\nu,k}}}
\frac{(h_{\nu,k}^2-x^2)\mathbf{H}_{\nu-1}(x)+2x\mathbf{H}_{\nu}(x)}{(h_{\nu,k}^2-x^2)\mathbf{H}_{\nu}(x)}\\
&=\frac{1}{2h_{\nu,k}}\lim_{x\rightarrow{h_{\nu,k}}}\frac{(h_{\nu,k}^2-x^2)\mathbf{H}_{\nu-1}''(x)-4x\mathbf{H}_{\nu-1}'(x)-2\mathbf{H}_{\nu-1}(x)+4\mathbf{H}_{\nu}'(x)+2x\mathbf{H}_{\nu}''(x)}
{(h_{\nu,k}^2-x^2)\mathbf{H}_{\nu}''(x)-4x\mathbf{H}_{\nu}'(x)-2\mathbf{H}_{\nu}(x)}\\
&=\frac{1}{2h_{\nu,k}}\left(\frac{-4h_{\nu,k}\mathbf{H}_{\nu-1}'(h_{\nu,k})-2\mathbf{H}_{\nu-1}(h_{\nu,k})+4\mathbf{H}_{\nu}'(h_{\nu,k})+2h_{\nu,k}\mathbf{H}_{\nu}''(h_{\nu,k})}
{-4h_{\nu,k}\mathbf{H}_{\nu}'(h_{\nu,k})}\right)\\
&=\frac{\nu-1}{2h_{\nu,k}^2}+\frac{h_{\nu,k}^{\nu-2}}{\sqrt{\pi}2^{\nu-1}\Gamma\left(\nu+\frac{1}{2}\right)\mathbf{H}_{\nu}'(h_{\nu,k})},
\end{align*}
which completes the proof of \eqref{formStruve}. Now, to prove \eqref{formStruve4} observe that by means of \eqref{formStruve} and the analogous of \eqref{vignat} for the zeros of Struve functions we have
\begin{align*}\sum_{n\geq 1,n\neq k}\frac{1}{h_{\nu,n}^4-h_{\nu,k}^4}&=\frac{1}{2h_{\nu,k}^2}\sum_{n\geq 1,n\neq k}\frac{1}{h_{\nu,n}^2-h_{\nu,k}^2}-\frac{1}{2h_{\nu,k}^2}\sum_{n\geq 1,n\neq k}\frac{1}{h_{\nu,n}^2+h_{\nu,k}^2}\\
&=\frac{1}{2h_{\nu,k}^2}\left(\frac{\nu+2}{2h_{\nu,k}^2}-\frac{h_{\nu,k}^{\nu-2}}{\sqrt{\pi}2^{\nu-1}\Gamma\left(\nu+\frac{1}{2}\right)\mathbf{H}_{\nu}'(h_{\nu,k})}\right)-
\frac{1}{2h_{\nu,k}^2}\left(\sum_{n\geq 1}\frac{1}{h_{\nu,n}^2+h_{\nu,k}^2}-\frac{1}{2h_{\nu,k}^2}\right),
\end{align*}
which completes the proof.
\end{proof}

\begin{proof}[\bf Proof of Theorem \ref{th3}] First we note that if $\nu=n+\frac12$ and $n\in\mathbb{N},$ then for all $z\in\mathbb{C}$
there holds
\begin{equation}\label{A0}\dfrac{2^{1-\nu}}{\Gamma(\nu)}z^\nu\,{\rm e}^{z}
K_\nu(z)=\prod_{k=1}^n\left(1-\dfrac{z}{z_{\nu,k}}\right),
\end{equation}
which by using logarithmic differentiation yields the next result
\begin{equation}\label{A5}
\dfrac{K_{\nu+1}(z)}{K_\nu(z)}=1+\dfrac{2\nu}{z}-\sum_{k=1}^n\dfrac{1}{z-z_{\nu,k}}.
\end{equation}
We note that the equality \eqref{A5} is already known \cite[Lemma 3.2]{YH} and
it is a corresponding Mittag--Leffler expansion for $K_\nu$. Now, let us notice that from \eqref{A5} it follows that
\begin{align*}
 \sum_{{k=1, k\neq
 j}}^n\dfrac{1}{z_{\nu,k}-z_{\nu,j}}&=\lim_{z\to
 z_{\nu,j}}\left(\dfrac{K_{\nu+1}(z)}{K_\nu(z)}-1-\dfrac{2\nu}{z}-\dfrac{1}{z_{\nu,j}-z}
 \right)\\
 &=\dfrac{1}{z_{\nu,j}}\lim_{z\to
 z_{\nu,j}}\dfrac{z(z_{\nu,j}-z)K_{\nu+1}(z)-(z+2\nu)(z_{\nu,j}-z)K_\nu(z)-zK_\nu(z)}{(z_{\nu,j}-z)\,K_\nu(z)}.
\end{align*}
On the other hand, by using the fact that $K_\nu$ is a particular solution of the modified Bessel
differential equation, i.e. there holds
\[z^2K_\nu''(z)+zK_\nu'(z)-(z^2+\nu^2)K_\nu(z)=0,\]
and in view of the relations
$$zK_\nu'(z)=\nu K_\nu(z)-zK_{\nu+1}(z),$$
\[zK_{\nu+1}'(z)+(\nu+1)K_{\nu+1}(z)=-zK_\nu(z),\]
we obtain the next identities for all $k\in\{1,2,\dots,n\}$
\begin{equation}\label{B1}
z_{\nu,k}K_\nu''(z_{\nu,k})+K_\nu'(z_{\nu,k})=0,
\end{equation}
\begin{equation}\label{B2}
K_\nu'(z_{\nu,k})=-K_{\nu+1}(z_{\nu,k})
\end{equation}
and
\begin{equation}\label{B3}
z_{\nu,k}K_{\nu+1}'(z_{\nu,k})+(\nu+1)K_{\nu+1}(z_{\nu,k})=0.
\end{equation}
Applying the Bernoulli-Hospital rule two times and in view of \eqref{B1}, \eqref{B2} and \eqref{B3} we
get
\begin{align*}
 \sum_{{k=1,k\neq
 j}}^n\dfrac{1}{z_{\nu,k}-z_{\nu,j}}&=
 \dfrac{-z_{\nu,j}K_\nu''(z_{\nu,j})-2z_{\nu,j}K_{\nu+1}'(z_{\nu,j})+2(z_{\nu,j}+2\nu-1)K_\nu'(z_{\nu,j})-
 2K_{\nu+1}(z_{\nu,j})}{-2z_{\nu,j}K_\nu'(z_{\nu,j})},\end{align*}
 which completes the proof of \eqref{eqk1}.

Now, we focus on the identities \eqref{eqk2} and \eqref{eqk4}. These can de deduced by using the next relations
\begin{align*}
\sum_{k=1,k\neq j}^n\frac{1}{z_{\nu,k}^2-z_{\nu,j}^2}&=\frac{1}{2z_{\nu,j}}\sum_{k=1,k\neq j}^n\frac{1}{z_{\nu,k}-z_{\nu,j}}-
\frac{1}{2z_{\nu,j}}\sum_{k=1,k\neq j}^n\frac{1}{z_{\nu,k}+z_{\nu,j}}\\
&=\frac{1}{2z_{\nu,j}}\frac{1-2\nu-2z_{\nu,j}}{2z_{\nu,j}}-\frac{1}{2z_{\nu,j}}\left(\sum_{k=1}^n\frac{1}{z_{\nu,k}+z_{\nu,j}}-\frac{1}{2z_{\nu,j}}\right)
\end{align*}
and
\begin{align*}
\sum_{k=1,k\neq j}^n\frac{1}{z_{\nu,k}^4-z_{\nu,j}^4}&=\frac{1}{2z_{\nu,j}^2}\sum_{k=1,k\neq j}^n\frac{1}{z_{\nu,k}^2-z_{\nu,j}^2}-
\frac{1}{2z_{\nu,j}^2}\sum_{k=1,k\neq j}^n\frac{1}{z_{\nu,k}^2+z_{\nu,j}^2}\\
&=\frac{1}{2z_{\nu,j}^2}\left(\frac{1-\nu-z_{\nu,j}}{2z_{\nu,j}^2}-\frac{1}{2z_{\nu,j}}\sum_{k=1}^n\frac{1}{z_{\nu,k}+z_{\nu,j}}\right)
-\frac{1}{2z_{\nu,j}^2}\left(\sum_{k=1}^n\frac{1}{z_{\nu,k}^2+z_{\nu,j}^2}-\frac{1}{2z_{\nu,j}^2}\right).\end{align*}

Alternatively, the identities \eqref{eqk2} and \eqref{eqk4} can be deduced as follows. First, let us notice that from \eqref{A0} it follows that
\begin{equation}\label{help1}
\prod_{k=1}^n\left(1-\dfrac{z^2}{z_{\nu,k}^2}\right)=\dfrac{2(-1)^\nu
z^{2\nu} K_\nu(z)K_\nu(-z)}{2^{2n}\Gamma(\nu)^2}.
\end{equation}
Logarithmic differentiation of the previous expression gives
\[-\sum_{k=1}^n\dfrac{2z}{z_{\nu,k}^2-z^2}=\dfrac{2\nu}{z}+\dfrac{K_{\nu}'(z)}{K_\nu(z)}-\dfrac{K_\nu'(-z)}{K_\nu(-z)},\]
which can be rewritten as follows
\[\sum_{k=1}^n\dfrac{1}{z_{\nu,k}^2-z^2}=\dfrac{K_{\nu+1}(z)}{2zK_\nu(z)}-\dfrac{K_{\nu+1}(-z)}{2zK_\nu(-z)}-\dfrac{2\nu}{z^2}.\]
This in turn implies that
\begin{equation}\label{B4}
 \sum_{{k=1,k\neq
 j}}^n\dfrac{1}{z_{\nu,k}^2-z_{\nu,j}^2}=\lim_{z\to
 z_{\nu,j}}\left(\dfrac{K_{\nu+1}(z)}{2zK_\nu(z)}-\dfrac{2\nu}{z^2}-\dfrac{1}{z_{\nu,j}^2-z^2}
 \right)-\dfrac{K_{\nu+1}(-z_{\nu,j})}{2z_{\nu,j}K_\nu(-z_{\nu,j})}=:\mathbf{K}_1-\mathbf{K_2}\,.
\end{equation}
By using the Bernoulli-Hospital rule,
two times, and then the relations \eqref{B1}, \eqref{B2} and
\eqref{B3}, we can conclude that
\begin{align}\label{B5}
\mathbf{K}_1&=\dfrac{1}{4z_{\nu,j}^3}\lim_{z\to
 z_{\nu,j}}\dfrac{z(z_{\nu,j}^2-z^2)K_{\nu+1}(z)-4\nu(z_{\nu,j}^2-z^2)
 K_\nu(z)-2z^2K_\nu(z)}{(z_{\nu,j}-z)\,K_\nu(z)}\\\nonumber
 &=\dfrac{-4z_{\nu,j}^2K_\nu''(z_{\nu,j})-4z_{\nu,j}^2K_{\nu+1}'(z_{\nu,j})-2(4z_{\nu,j}-8\nu
 z_{\nu,j})K_\nu'(z_{\nu,j})-6z_{\nu,j}K_{\nu+1}(z_{\nu,j})}{-8z_{\nu,j}^3K_\nu'(z_{\nu,j})}\\\nonumber
 &=\dfrac{1-3\nu}{2z_{\nu,j}^2}\,.
\end{align}
Using the substitution $z\mapsto-z_{\nu,j}$ from \eqref{A5} it
follows that
\begin{equation}\label{B6}
\mathbf{K}_2=\dfrac{1}{2z_{\nu,j}}-\dfrac{\nu}{z_{\nu,j}^2}+\dfrac{1}{2z_{\nu,j}}\sum_{k=1}^n\dfrac{1}{z_{\nu,k}+z_{\nu,j}}.
\end{equation}
Combination of \eqref{B4}, \eqref{B5} and \eqref{B6} yields the
proof of \eqref{eqk2}.

Similarly, let us notice that from \eqref{help1}, using the
substitution $z\mapsto {\rm i}\,z$ it follows that
\begin{equation*}
\prod_{k=1}^n\left(1+\dfrac{z^2}{z_{\nu,k}^2}\right)=\dfrac{2
z^{2\nu} K_\nu({\rm i}\,z)K_\nu(-{\rm
i}\,z)}{2^{2n}\Gamma(\nu)^2},
\end{equation*}
which in turn implies that
\begin{equation*}
\prod_{k=1}^n\left(1-\dfrac{z^4}{z_{\nu,k}^4}\right)=\dfrac{4(-1)^\nu
z^{4\nu} K_\nu(z)\,K_\nu(-z)\,K_\nu({\rm i}\,z)\,K_\nu(-{\rm
i}\,z)}{2^{4n}\Gamma(\nu)^4}.
\end{equation*}
Logarithmic differentiation of the previous expression gives
\[-\sum_{k=1}^n\dfrac{4z^3}{z_{\nu,k}^4-z^4}=\dfrac{4\nu}{z}+\dfrac{K_{\nu}'(z)}{K_\nu(z)}-\dfrac{K_\nu'(-z)}{K_\nu(-z)}+\dfrac{{\rm i}\,K_{\nu}'({\rm i}\,z)}{K_\nu({\rm i}\,z)}-
\dfrac{{\rm i}\,K_\nu'(-{\rm i}\,z)}{K_\nu(-{\rm i}\,z)},\] which can be rewritten as
\begin{equation*}
\sum_{k=1}^n\dfrac{4z^3}{z_{\nu,k}^4-z^4}=-\dfrac{8\nu}{z}+\dfrac{K_{\nu+1}(z)}{K_\nu(z)}-\dfrac{K_{\nu+1}(-z)}{K_\nu(-z)}+\dfrac{{\rm
i}\,K_{\nu+1}({\rm i}\,z)}{K_\nu({\rm i}\,z)}- \dfrac{{\rm
i}\,K_{\nu+1}(-{\rm i}\,z)}{K_\nu(-{\rm i}\,z)}\,.
\end{equation*}
Consequently we have
\begin{align}\label{B44}
\sum_{k=1, k\neq j}^n\dfrac{1}{z_{\nu,k}^4-z_{\nu,j}^4}&=\dfrac{1}{4z_{\nu,j}^3}\lim_{z\to
 z_{\nu,j}}\left(\dfrac{K_{\nu+1}(z)}{K_{\nu}(z)}-\dfrac{8\nu}{z}-\dfrac{4z^3}{z_{\nu,j}^4-z^4}
 \right)\\\nonumber&\qquad+\dfrac{1}{4z_{\nu,j}^3}\left( -\dfrac{K_{\nu+1}(-z_{\nu,j})}{K_\nu(-z_{\nu,j})}+\dfrac{{\rm i}\,K_{\nu+1}({\rm i}\,z_{\nu,j})}{K_\nu({\rm i}\,z_{\nu,j})}-
\dfrac{{\rm i}\,K_{\nu+1}(-{\rm i}\,z_{\nu,j})}{K_\nu(-{\rm
i}\,z_{\nu,j})}\right)=:\mathbf{K}_3+\mathbf{K}_4\,.
\end{align}
By using again the Bernoulli-Hospital rule two times, and then the relations
\eqref{B1}, \eqref{B2} and \eqref{B3}, we obtain that
\begin{align}\label{B55}
\mathbf{K}_3&=\dfrac{1}{16z_{\nu,j}^7}\lim_{z\to
 z_{\nu,j}}\dfrac{z(z_{\nu,j}^4-z^4)K_{\nu+1}(z)-8\nu(z_{\nu,j}^4-z^4)
 K_\nu(z)-4z^4K_\nu(z)}{(z_{\nu,j}-z)\,K_\nu(z)}\\\nonumber
 &=\dfrac{-4z_{\nu,j}^4K_\nu''(z_{\nu,j})-8z_{\nu,j}^4K_{\nu+1}'(z_{\nu,j})+32z_{\nu,j}^3(2\nu-1)K_\nu'(z_{\nu,j})-20z_{\nu,j}^3K_{\nu+1}(z_{\nu,j})}{-32z_{\nu,j}^7K_\nu'(z_{\nu,j})}\\\nonumber
 &=\dfrac{2-7\nu}{4z_{\nu,j}^4}\,.
\end{align}
Finally, from \eqref{A5} it holds
\begin{align}\label{B66}
\mathbf{K}_4&=\dfrac{1}{4z_{\nu,j}^3}\left(-1+\dfrac{6\nu}{z_{\nu,j}}-
 \sum_{k=1}^n\left(\dfrac{1}{z_{\nu,j}+z_{\nu,k}}+\dfrac{{\rm i}}{{\rm i}\,z_{\nu,j}-z_{\nu,k}}+\dfrac{{\rm i}}{{\rm
 i}\,z_{\nu,j}+z_{\nu,k}}\right)\right)\\\nonumber &=\dfrac{1}{z_{\nu,j}^3}\left(-1+\dfrac{6\nu}{z_{\nu,j}}-
 \sum_{k=1}^n\left(\dfrac{1}{z_{\nu,j}+z_{\nu,k}}+\dfrac{2z_{\nu,j}}{z_{\nu,j}^2+z_{\nu,k}^2}\right)\right)\,.
\end{align}
From \eqref{B44}, \eqref{B55} and \eqref{B66} the desired formula \eqref{eqk4} immediately follows.
\end{proof}

\subsection*{Acknowledgement} The first author is grateful to Christophe Vignat for providing an alternative proof of \eqref{formnew} as well as for the useful discussions during the process of writing this paper.


\begin{thebibliography}{}

\bibitem{bps}
\textsc{\'A. Baricz, S. Ponnusamy, S. Singh}, Tur\'an type inequalities for Struve functions, {\em Proc. Amer. Math. Soc.} (submitted).

\bibitem{calogero}
\textsc{F. Calogero}, On the zeros of Bessel functions, {\em Lett. Nuovo Cimento} 20(7) (1977) 254--256.

\bibitem{FS}
\textsc{E.M. Ferreira, J. Sesma}, Zeros of the Macdonald function of complex order, {\em J. Comput. Appl. Math.} 211 (2008) 223--231.

\bibitem{YH}
\textsc{Y. Hamana}, The expected volume and surface area of the Wiener sausage in odd dimensions,
{\em Osaka J. Math.} 49 (2012) 853--868.

\bibitem{HMS} \textsc{Y. Hamana, H. Matsumoto, T. Shirai}, On the zeros of the Macdonald functions,
arXiv: 1302.5154v1.

\bibitem{ismail}
\textsc{M.E.H. Ismail, M.E. Muldoon}, Certain monotonicity properties of Bessel functions, {\em J. Math. Anal. Appl.} 118 (1986) 145--150.

\bibitem{nist} \textsc{F.W.J. Olver, D.W. Lozier, R.F. Boisvert, C.W. Clark} (Eds.), {\em NIST Handbook of Mathematical Functions}, Cambridge Univ. Press, Cambridge, 2010.

\bibitem{Watson}
\textsc{G.N. Watson}, {\em A Treatise on the Theory of Bessel
Functions}, Cambridge Univ. Press, Cambridge, 1944.
\end{thebibliography}
\end{document}